\documentclass{lmcs}
\pdfoutput=1

\usepackage{lastpage}
\lmcsdoi{18}{1}{7}
\lmcsheading{}{\pageref{LastPage}}{}{}%
{Jan.~30,~2018}{Jan.~12,~2022}{}

\usepackage[utf8]{inputenc}


\usepackage{amssymb}

\def\sub{\subseteq}

\def\cov{\lhd}

\def\Pos{\mathrm{\,Pos}}
\def\P{\mathcal{P}_{\omega_1}}
\def\N{\mathbb{N}}

\begin{document}

\title[\texorpdfstring{$\sigma$}{sigma}-locales in Formal Topology]{\texorpdfstring{$\sigma$}{sigma}-locales in Formal Topology\rsuper*}

\titlecomment{{\lsuper*}This project has received funding from the European Union’s Horizon 2020 research and innovation programme under the Marie Skłodowska-Curie grant agreement No 731143}

\author[F.~Ciraulo]{Francesco Ciraulo}
\address{Department of Mathematics, University of Padua, Via Trieste 63, 35121 Padova (Italy)}
\email{ciraulo@math.unipd.it}
\urladdr{\url{www.math.unipd.it/~ciraulo}}

\begin{abstract}
\noindent A $\sigma$-frame is a poset with countable joins and finite meets in which binary meets distribute over countable joins.
The aim of this paper is to show that $\sigma$-frames, actually $\sigma$-locales, can be seen as a branch of Formal Topology, that is, intuitionistic and predicative point-free topology. Every $\sigma$-frame $L$ is the lattice of Lindel\"of elements (those for which each of their covers admits a countable subcover) of a formal topology of a specific kind which, in its turn, is a presentation of the free frame over $L$. We then give a constructive characterization of the smallest (strongly) dense $\sigma$-sublocale of a given $\sigma$-locale, thus providing a ``$\sigma$-version'' of a Boolean locale. Our development depends on the axiom of countable choice.
\end{abstract}

\maketitle

\section*{Introduction}

It is well known that the set $B(H) = \{x\in H\ |\ x=--x\}$ of stable elements of a complete Heyting algebra $H$ is a complete Boolean algebra. Actually $B(H)$ is a quotient of $H$ in the category of frames. From the point of view of the category of locales, this means that every locale $L$ contains a Boolean sublocale $B(L)$, which can be characterized as the smallest dense sublocale of $L$.

Sambin introduced the notion of an \emph{overlap algebra} (see~\cite{regolari, formal_regular}) as a ``positive'' alternative to that of a complete Boolean algebra. One of the main advantages of his approach is that powersets are examples of overlap algebras (in fact they are precisely the atomic ones), although they are not Boolean, constructively.

It has recently turned out~\cite{almost_discrete} (see also~\cite{contente}) that overlap algebras can be understood as the smallest \emph{strongly} dense sublocales (in the sense of~\cite{strong_density}) of overt locales. The same statement can be given a predicative interpretation by considering a formal topology $(S,\cov,\Pos)$ in place of an overt locale $L$.

The notion of a $\sigma$-locale is a natural generalization of that of a locale: the underlying lattice is a $\sigma$-frame, rather than a frame, that is, it is required to have just countable, rather than arbitrary, joins.  As shown in~\cite{simpson}, $\sigma$-locales play an important role in the point-free approach to measure theory and probability.

The construction of $B(L)$ from $L$ can be mimicked in the case of $\sigma$-locales~\cite{madden}. In that case, $B(L)$ is still the smallest dense $\sigma$-sublocale of $L$; however, it is not Boolean any longer, in general. The $\sigma$-frames of the form $B(L)$ are called \emph{d-reduced} (``d'' for ``dense'') in~\cite{madden}.

One of our aims is to give a positive account of d-reduced $\sigma$-locales. In order to obtain this, we work with $\sigma$-locales which are \emph{overt} (in a suitable sense). The positivity predicate $\Pos$ of an overt $\sigma$-locale $L$ is then used to define a positive version of the \emph{codense} congruence relation on $L$~\cite{madden}, which corresponds to the smallest dense $\sigma$-sublocale $B(L)$ of $L$. Actually, because of the positive nature of our definition, the notion of density involved here is intuitionistically stronger than the usual one (in accordance with what happens in the case of locales, as mentioned above).

Our arguments are always intuitionistically valid and predicative; but we need the axiom of countable choice. In fact, our results could be formalized in the extensional level of the so-called Minimalist Foundations~\cite{maietti-sambin,maietti} augmented with countable choice. In such a foundational framework, Formal Topology is the ``native'' way to develop point-free topology.

The paper is organized as follows. In Section 1 we recall some constructive results about the notion of a countable set. Section 2 deals with $\sigma$-frames and $\sigma$-locales within the framework of Formal Topology. Finally, Section 3 presents the construction of the smallest strongly dense $\sigma$-sublocale of an overt $\sigma$-locale.

\section{A constructive look at countable sets}

By a countable set we intuitively mean a set $S$ which is either (empty or) finite or countably infinite, that is, in bijection with the set $\N $ of natural numbers. Within usual foundations (such as $ZF$ with countable choice $AC_\omega$), this is equivalent to saying that a set is either empty or enumerable in the sense that there exists an onto map $\N \twoheadrightarrow S$. This case distinction looks inappropriate for a good constructive definition. Following a quite established tradition (see for instance~\cite{bauer}), we give the following (seemingly tricky but, as we will see, definitely convenient) definition.

\begin{defi}
A set $S$ is \emph{countable} if there exists a surjection $\alpha:\N \twoheadrightarrow S+\{\bot\}$. Equivalently, $S$ is \emph{countable} if there exists a map $\alpha:\N \to S+\{\bot\}$ such that $S\sub\alpha[\N ]$.
\end{defi}

\noindent Here $S+\{\bot\}$ is the disjoint union (or sum) of $S$ and $\{\bot\}$, and $\alpha[\N ]$ is the image of $\N $ along $\alpha$. For the sake of notational simplicity, we do not distinguish between an element of $S$ and its copy inside $S+\{\bot\}$; otherwise, we should have written $\forall a\in S.\exists n\in\N .\alpha(n)=i(a)$, where $i$ is the canonical injection of $S$ into $S+\{\bot\}$, instead of the more readable $S\sub\alpha[\N ]$.

\begin{rem} The term ``Constructive Mathematics'' refers, as it is well known, to a variety of foundational approaches ranging from intuitionistic type theories to constructive set theories, from topos valid mathematics to constructive mathematics à la Bishop. Pragmatically, we shall try to keep ourselves as neutral as possible with respect to the different foundational choices and hence to provide definitions and proofs in such a way that they remain valid and meaningful within virtually any foundations. In this paper this would be possible only up to a certain extent, as we shall need the Axiom of Countable Choice $AC_\omega$
\[
\forall n\in\N .\exists x\in X.R(n,x)\Longrightarrow\exists\alpha:\N \to X.\forall n\in\N .R(n,\alpha(n))
\]
for every set $X$ and every relation $R\sub\N\times X$.
Apart from this, our position automatically forces us to abandon the so-called Law of Excluded Middle (LEM), the full Powerset Axiom (PA), and the full Axiom of Choice (AC). Also, we will have to consider general collections which cannot be assumed to be sets, an important example being powersets. Of course, only a formal theory can clarify what is precisely meant by a set. Here we just need the existence of the set of natural numbers together with its initial (finite) segments, and we need the class of sets to be closed under disjoint union, exponentiation and quotients.
If asked for a concrete theory to be assumed as a foundations, we would suggest the (extensional level of the) Minimalist Foundation~\cite{maietti, maietti-sambin} augmented with $AC_\omega$.
\end{rem}

A subset $D\sub S$ is {\bf detachable} if there exists an operation $\chi_D:S\to 2$ (the characteristic function of $D$), where $2$ is the set $\{0,1\}$ of Boolean values, such that $x\in D\Leftrightarrow \chi_D(x)=1$.

\begin{prop}
A set $S$ is countable if and only if there exists a surjective map $D\twoheadrightarrow S$ with $D$ a detachable subset of $\N $.
\end{prop}
\begin{proof}
Assume that $S$ is countable and that $\alpha:\N \to S+\{\bot\}$ is the ``evidence'' of that (as required by the definition). Define $D$ = $\{n\in\N \ |\ \alpha(n)\in S\}$, which is detachable (because one can decide to which part of a disjoint union an element belongs). The restriction of $\alpha$ to $D$ is a surjection onto $S$.

Vice versa, given $g:D\twoheadrightarrow S$, define $\alpha(n)$ as either $g(n)$ or $\bot$ according to whether $n$ belong to $D$ or not (that is, according to whether $\chi_D(n)$ is 1 or 0).
\end{proof}

Such a characterization has in fact been taken as a definition in~\cite{varieties}. Classically, of course, every set of natural numbers is detachable and so the previous proposition says just that $S$ is countable if and only if its cardinality is not greater than $\aleph_0$.

\begin{lem}\label{lemmaDetachableCountable} Every detachable subset of a countable set is countable.
\end{lem}
\begin{proof}
Let $\alpha:\N \to S+\{\bot\}$ be such that $S\sub\alpha[\N ]$, and let $X$ be a detachable subset of $S$. Define $\beta:\N \to S+\{\bot\}$ as follows: put $\beta(n)=\alpha(n)$ if $\alpha(n)\in X$, and put $\beta(n)=\bot$ otherwise. Clearly $S\cap\beta[\N]=X$.
\end{proof}

Note that there cannot be a general way to decide whether a countable set is inhabited or not: this would imply the Limited Principle of Omniscience LPO (see Remark~\ref{remark_Sigma} below). For a Brouwerian counterexample, consider the set of even numbers greater than 4 which are not the sum of two odd primes: it is detachable, hence countable by Lemma~\ref{lemmaDetachableCountable}, but we still do not know if it is empty.

\subsection{The set of countable subsets}

Given a set $S$, a subset $X\sub S$ is a countable set if and only if there exists $\alpha:\N \to S+\{\bot\}$ such that $X=S\cap\alpha[\N ]$.
We write $\P(S)$ for the collection of all countable subsets of $S$. Clearly we have
\[\P(S)\quad\cong\quad(S+\{\bot\})^\N /\sim\]
where $\alpha\sim\beta$ means $S\cap\alpha[\N ]\ =\ S\cap\beta[\N ]$. Hence $\P(S)$ is a set (it is a quotient of a set).\footnote{On the contrary, we do not assume that $\mathcal{P}(S)$, the collection of all subsets of $S$, is a set.} Note that the set of (Kuratowski-)finite subsets of $S$ (see, for instance,~\cite{finiteness}) can be identified with a subset of $\P(S)$. 

A set $S$ has a {\bf decidable equality} if the diagonal $\{(a,a)\ |\ a\in S\}$ is a detachable subset of $S\times S$.

\begin{prop}\label{prop1} For every set $S$, $\P(S)$ is closed under countable unions. And if equality in $S$ is decidable, then $\P(S)$ is closed under binary intersections.
\end{prop}
\begin{proof}
Let $\{X_i\ |\ i\in I\}$ be a countable family of countable subsets of $S$. So there exists $\alpha:\N \to I+\{\bot\}$ such that $I\sub\alpha[\N ]$. For each $i\in I$, we choose (by $AC_\omega$) a map $\beta_i:\N \to S+\{\bot\}$ such that $X_i=S\cap\beta_i[\N ]$. We want to check that $\bigcup_{i\in I}X_i$ is countable. Indeed, it is enumerated by the map $\gamma:\N \to S+\{\bot\}$ defined as follows. First, by means of a suitable (recursive) pairing function, we identify $\N$ with $\N \times\N$. Second, we define the image of the pair $(n,m)$ to be
$\beta_{\alpha(n)}(m)$ if $\alpha(n)\in I$, and $\bot$ otherwise.

As for the second part of the statement, note that $\P(S)$ is closed under binary intersections if and only if $\{a\}\cap\{b\}$ is countable for every $a,b\in S$. Indeed, given two countable subsets $X_i=\{a_{i,n}\ |\ n\in D_i\}$, $i=1,2$, their intersection $X_1\cap X_2$ can be written as the countable union $\bigcup_{n\in D_1}\bigcup_{m\in D_2}(\{a_{1,n}\}\cap\{a_{2,m}\})$. Now if equality is decidable, then $\{a\}\cap\{b\}$ is either empty or a singleton, and hence it is countable.
\end{proof}

The special case $\P(1)$, where $1=\{0\}$, is sometimes written $\Sigma$; it is a subset of the collection $\Omega$ of all truth values (that is, the collection $\mathcal{P}(1)$ of all subsets of $1$). In fact, an element of $\Sigma$ can be identified with (the truth value of) a proposition of the form $\exists n.[\alpha(n)=0]$, for some $\alpha:\N \to 1+\{\bot\}$. Equivalently, an element of $\Sigma$ can be thought of as (the truth value of) the proposition ``$D$ is inhabited'', for some detachable $D\sub\mathbb{N}$. So $\Sigma$ is precisely what is known as the {\bf Rosolini dominance}~\cite{rosolini}; it is the set of ``open'' (or ``semi-decidable'') truth values in Synthetic Topology~\cite{bauer}.

\begin{rem}\label{remark_Sigma} Classically, of course, $\Sigma$ = 2. Constructively, $\Sigma$ = 2 is equivalent to requiring that every $p\in\Sigma$ is either inhabited or empty; and this is equivalent to LPO, that is, the assertion $\forall f:\N\to 2.(\exists n\in\mathbb{N}.f(n)=0 \vee\forall n\in\mathbb{N}.f(n)=1)$. Indeed, by interpreting 2 as $1+\bot$, LPO becomes $\forall p\in\Sigma.(p\vee\neg p)$.
\end{rem}

The second part of Proposition~\ref{prop1} can be strengthened, as we are now going to show.
We say that a set $S$ has a {\bf semi-decidable equality} if there exists an operation $\psi:S\times S\to\Sigma$ such that $a=b\Longleftrightarrow \psi(a,b)=1$.

\begin{prop}\label{prop2}For every set $S$ the following are equivalent:
\begin{enumerate}
\item $\P(S)$ is closed under binary intersections;
\item equality in $S$ is semi-decidable.
\end{enumerate}
\end{prop}
\begin{proof}
If $\P(S)$ is closed under binary intersections, then we can define a map $\psi:S\times S\to\Sigma$ by putting $\varphi(a,b)=\{x\in 1\ |\ a=b\}$. We claim that this is a countable subset of 1. By assumption there exists $D\twoheadrightarrow\{a\}\cap\{b\}$, with $D\sub\N$ detachable, which we can compose with the obvious map from $\{a\}\cap\{b\}$ onto $\{x\in 1\ |\ a=b\}$.\footnote{Given any set $S$, there exists precisely one map from $S$ to the terminal set $1=\{0\}$. Such a map factorizes via its image, which is just $\{x\in 1\ |\ S\textrm{ is inhabited}\}$.}

Vice versa, given any $a,b\in S$, we have $\psi(a,b)\in\Sigma$ and so there exists $f:\N\to 2$ such that $\psi(a,b)$ is the truth value of $\exists n\in\N.f(n)=0$. We can use $f$ to define a map $\alpha:\N\to S+\{\bot\}$ as follows: we put $\alpha(n)=a(=b)$ if $f(n)=0$ and $\alpha(n)=\bot$ otherwise. Clearly $S\cap\alpha[\N]$ = $\{a\}\cap\{b\}$.
\end{proof}

Finally, note that $\P(S)$ has a largest element if and only if $S$ itself is countable.

\begin{rem}\label{remark_inhabited}
For future reference, note that the statement ``$W$ is inhabited'' for $W\in\P(S)$ can be seen as an element of $\Sigma$. Indeed, if $W=\{a_i\ |\ i\in D\}$ for some detachable $D\sub\N$, then ``$W$ is inhabited'' is equivalent to ``$D$ is inhabited''.
\end{rem}

\section{\texorpdfstring{$\sigma$}{sigma}-frames in Formal Topology}

A {\bf suplattice} is a partially ordered collection $(P,\leq)$ with all set-indexed joins (hence a bottom element). Here $P$ need not be a set. A {\bf base} for $P$ is a subset $S\sub P$ such that, for every $p\in P$, (i) $\{a\in S\ |\ a\leq p\}$ is a set and (ii) $p=\bigvee\{a\in S\ |\ a\leq p\}$. In that case, $(P,\leq,S)$ is called a {\bf set-based} suplattice.\footnote{A set-based suplattice has all set-indexed meets too (hence a top element) because $\bigwedge_{i\in I}p_i$ is just $\bigvee\{a\in S\ |\ a\leq p_i\textrm{ for all }i\in I\}$.} All the information about $(P,\leq,S)$ can be encoded as a pair $(S,\cov)$ where $a\cov U$ is $a\leq\bigvee U$, for $a\in S$ and $U\sub S$. The structure $(S,\cov)$ is called a {\bf basic cover} and it is characterized abstractly by the following two properties:
\begin{enumerate}
\item $a\cov U$ whenever $a\in U$, and
\item if $a\cov U$ and $u\cov V$ for all $u\in U$, then $a\cov V$.
\end{enumerate}
Given $(S,\cov)$, the suplattice $P$ can be recovered (up to isomorphism) as a quotient $\mathcal{P}(S)/=_\cov$, where $U=_\cov V$ means $\forall a\in S.(a\cov U\Leftrightarrow a\cov V)$. A join $\bigvee_{i\in I}[U_i]$ in $\mathcal{P}(S)/=_\cov$ is computed as $[\bigcup_{i\in I}U_i]$; in particular, $[U]\leq[V]$ if and only if $\forall a\in U.a\cov V$.

A {\bf frame} is a suplattice with finite meets in which binary meets distribute over (set-indexed) joins. Set-based frames correspond to a special class of basic covers called {\bf formal covers}. Actually there are a number of different ways to explicitly define the notion of a formal cover~\cite{convergence}; in all cases, of course, the resulting category is (dually) equivalent to that of set-based frames. Here we prefer the following definition which corresponds to assuming the base $S$ to be closed under finite meets (a property that can always be assumed for every set-based frame without loss of generality).

\begin{defi}
A \emph{formal cover} is given by a basic cover $(S,\cov)$ together with an inf-semilattice structure $(S,\wedge,1)$ such that
\begin{enumerate}
\item $a\cov\{1\}$
\item $a\cov U$ $\Longrightarrow$ $(a\wedge b)\cov\{u\wedge b\ |\ u\in U\}$
\end{enumerate}
for all $a,b\in S$ and $U\sub S$.\footnote{A notable consequence of this definition is that $a\leq b$ $\Rightarrow$ $a\cov\{b\}$; indeed $a\leq b$ means that $a$ equals $a\wedge b$, and $a\wedge b\cov\{b\}$ follows from $a\cov\{1\}$ (since $1\wedge b$ equals $b$).}
\end{defi}

Given a formal cover, binary meets in the corresponding frame $\mathcal{P}(S)/=_\cov$ are computed as $[U]\wedge[V]$ = $[\{u\wedge v\ |\ u\in U,\ v\in V\}]$ in terms of the meet operation of $S$; moreover $[\{1\}]=[S]$ is the top element of the frame.

A {\bf $\sigma$-frame} is a partial order with countable joins and finite meets, in which binary meets distribute over countable joins. In this paper, we restrict our attention to $\sigma$-frames whose carriers are sets. For instance, $\P(S)$ is a $\sigma$-frame if $S$ has semi-decidable equality (Propositions~\ref{prop1} and~\ref{prop2}) and if, at the same time, $S$ is countable. A homomorphism of $\sigma$-frames is a map which preserves countable joins and finite meets.

It is easy to see that $\Sigma=\P(1)$ is initial in the category of $\sigma$-frames (actually, $\Sigma$ is the free $\sigma$-frame on no generators).

\subsection{\texorpdfstring{$\sigma$}{sigma}-coherent formal topologies}

Formal covers are a powerful tool, for instance when it comes to constructing the free frame over a given $\sigma$-frame. This is done in this section.

Let $L$ be a $\sigma$-frame. For $a\in L$ and $U\sub L$, let us put
\begin{equation}\label{eq.cover_L}
a\cov_L U\stackrel{def}{\Longleftrightarrow}a\leq\bigvee W\textrm{ for some countable subset }W\sub U
\end{equation}
(note that $a\cov_L\{b\}$ holds if and only if $a\leq b$ in $L$).

\begin{prop}
For $L$ a $\sigma$-frame, the pair $(L,\cov_L)$ is a formal cover.
\end{prop}
\begin{proof}
Checking that $(L,\cov_L)$ is a formal cover is quite straightforward. Only one point requires the axiom of countable choice, namely showing that if $a\cov_L U$ and $u\cov_L V$ for all $u\in U$, then $a\cov_L V$. Indeed we have $a\leq\bigvee W$ for some countable $W\sub U$ and, for each $u\in W$, we also have $u\leq\bigvee W_u$ for some countable $W_u\sub V$. So $a\leq\bigvee W\leq\bigvee_{u\in W}\bigvee W_u=\bigvee\bigcup_{u\in W}W_u$. Since $\bigcup_{u\in W}W_u$ is a countable subset of $V$, we can conclude that $a\cov_L V$.
\end{proof}

The frame presented by $(L,\cov_L)$ is called the {\bf frame envelope} of $L$ in~\cite{envelope}.
Some results about the frame envelope look quite elegant in the language of formal topology. For instance, the statement that if $L$ is compact as a $\sigma$-frame, then its envelope is compact as a frame~\cite{envelope} can be proved as follows.\footnote{By directly interpreting the notion proposed in~\cite{envelope} in our foundational framework, we say that a $\sigma$-frame $L$ is compact if, for any countable $W\subseteq L$, $1=\bigvee W$ implies that $1=\bigvee K$ for some Kuratowski-finite $K\subseteq W$.} Assume $1\cov_L U$, that is, $1\leq\bigvee W$ for some countable $W\sub U$. By compactness of $L$, there is a (Kuratowski-)finite $K\sub W\sub U$ such that $1\leq\bigvee K$. So $1\cov_L K$ and hence $(L,\cov_L)$ is compact. 

The frame envelope $P$ = $\mathcal{P}(L)/=_{\cov_L}$ is the free frame over $L$ as a $\sigma$-frame. In general, a frame $P$ is the\footnote{As usual, this must be understood up to isomorphism.} {\bf free frame} over the $\sigma$-frame $L$ if there exists a $\sigma$-frame homomorphism $m:L\to P$ such that for every $\sigma$-frame homomorphism $f:L\to Q$ with $Q$ a frame there exists a unique frame homomorphism $h:P\to Q$ with $h\circ m$ = $f$. In the case of the frame envelope $P$ of $L$, $m$ and $h$ are defined as $m(a)=[a]$ and $h([U])$ = $\bigvee_{a\in U}f(a)$.

Free frames over the category of $\sigma$-frames, as constructed in the previous proposition, can be characterized explicitly as follows.

Given a frame $P$, say that $a\in P$ is {\bf Lindel\"of}~\cite{envelope} if
\begin{equation}
a\leq\bigvee X\Longrightarrow a\leq\bigvee W\textrm{ for some countable }W\sub X
\end{equation}
for all $X\sub P$. Lindel\"of elements are closed under countable joins (by $AC_\omega$), but not under finite meets, in general.

A frame $P$ is called {\bf $\sigma$-coherent}~\cite{madden} if
\begin{enumerate}
\item its Lindel\"of elements form a set\footnote{The requirement that the Lindel\"of elements form a set seems necessary in our framework.} and they are closed under finite meets (hence they form a $\sigma$-frame), and
\item every element of $P$ is a (not necessarily countable) join of Lindel\"of elements.
\end{enumerate}
For instance, if $S$ is countable with semi-decidable equality, then $\mathcal{P}(S)$ is $\sigma$-coherent and its $\sigma$-frame of Lindel\"of elements is just $\P(S)$. In particular, $\Omega=\mathcal{P}(1)$ is $\sigma$-coherent and $\Sigma$ is its $\sigma$-frame of Lindel\"of elements.

The set $L$ of Lindel\"of elements of a $\sigma$-coherent frame $P$ is a base for $P$. Therefore $P$ can be presented as a formal cover $(L,\cov)$ where, as usual, $a\cov U$ means $a\leq\bigvee U$. However, since $a\in L$ is Lindel\"of, $a\cov U$ happens precisely when $a\leq\bigvee W$ for some countable $W\sub U$. In other words, $\cov$ is just $\cov_L$ as defined in the previous proposition. This immediately gives the following result.

\begin{prop} For a frame $P$, the following are equivalent:
\begin{enumerate}
\item $P$ is $\sigma$-coherent;
\item $P$ is the frame envelope $(L,\cov_L)$ of some $\sigma$-frame $L$, in which case $L$ is (isomorphic to) the $\sigma$-frame of Lindel\"of elements of $P$;
\item $P$ is the free frame over some $\sigma$-frame $L$.\qed%
\end{enumerate}
\end{prop}

\noindent
Formal covers of the form $(L,\cov_L)$ are characterized, up to isomorphism, as those formal covers which satisfy the following equation~\eqref{eq-sigma_cover}. This fact follows immediately from the previous discussion and the next proposition.

\begin{defi}
A formal cover $(S,\cov)$ is called a {\bf $\sigma$-cover} if
\begin{equation}\label{eq-sigma_cover}
a\cov U\Longrightarrow a\cov W\textrm{ for some countable subset }W\sub U
\end{equation}
for every $a\in S$ and $U\sub S$.
\end{defi}

\begin{prop}
The frame presented by a $\sigma$-cover $(S,\cov)$ is $\sigma$-coherent. The corresponding $\sigma$-frame of Lindel\"of elements is $\P(S)/=_\cov$.
\end{prop}
\begin{proof}
Condition~\eqref{eq-sigma_cover} says that $[\{a\}]$ is Lindel\"of for every $a\in S$ (in particular the top element $[\{1\}]$ is Lindel\"of). Therefore, $[W]$ is Lindel\"of for every countable $W\subseteq S$ (this requires $AC_\omega$). Moreover, an element $[U]$ is Lindel\"of (if and) only if $[U]$ = $[W]$ for some countable $W$. Indeed, from $[U]$ $\leq$ $\bigvee_{a\in U}[\{a\}]$ one gets $[U]$ $\leq$ $\bigvee_{a\in W}[\{a\}]$ = $[W]$ for some countable $W\sub U$. So the collection of Lindel\"of elements can be identified with the set $\P(S)/=_\cov$.\footnote{There is a subtlety in this proof. In the attempt to identify the collection of Lindel\"of elements of a $\sigma$-cover $(S,\cov)$ with the set $\P(S)/=_\cov$, we can consider the inclusion of the latter into the former, which turns out to be onto by the definition of a $\sigma$-cover. However, if we want to be able to construct the inverse mapping from the Lindel\"of elements to $\P(S)/=_\cov$, it seems necessary (if we want to avoid the axiom of choice) to strengthen the definition of a $\sigma$-cover by requiring an explicit operation which computes a countable subset $W$ for each $a\in S$ and $U\subseteq S$ with $a\cov U$.} If $W_1$ and $W_2$ are countable, then $\{w_1\wedge w_2\ |\ w_1\in W_1,\ w_2\in W_2\}$ is countable too (by pairing), so that Lindel\"of elements are closed under binary meets.\footnote{Note that we are not assuming that $S$ is countable with semi-decidable equality here, so $\P(S)$ need not be closed under finite intersections.}
\end{proof}

\subsection{\texorpdfstring{$\sigma$}{sigma}-locales as formal topologies}

The previous results show that a $\sigma$-coherent frame is essentially the same thing as a $\sigma$-frame (namely the $\sigma$-frame of its Lindel\"of elements). This suggests to present the category of $\sigma$-frames as a (non full) subcategory of the category of frames. In order to do that, one has to consider only those frame homomorphisms which preserve Lindel\"of elements (by freeness, each $\sigma$-frame homomorphism is the restriction of a unique frame homomorphism between the corresponding envelopes).

Here we prefer to work with the category of {$\sigma$-locales}, the opposite of the category of $\sigma$-frames. (However, we make no notational distinction between a $\sigma$-locale and its corresponding $\sigma$-frame.)

If $(S_1,\cov_1)$ and $(S_2,\cov_2)$ are basic covers, then a suplattice homomorphism $h$ from $\mathcal{P}(S_2)/=_{\cov_2}$ to $\mathcal{P}(S_1)/=_{\cov_1}$ (note the direction) can be presented by means of a binary relation $r\sub S_1\times S_2$ such that $U=_{\cov_2} V$ $\Rightarrow$ $r^{-1}U=_{\cov_1}r^{-1}V$ for all $U,V\sub S_2$, where $r^{-1}Y=\{x\in S_1\ |\ r(x,y)\textrm{ for some }y\in Y\}$. Indeed, every such $r$ induces a homomorphism $h$ given by $h([Y])$ = $[r^{-1}Y]$; vice versa, given $h$ it is sufficient to put $r(x,y)$ iff $x\cov_1 h([\{y\}])$. Note that several relations may be used to define the same homomorphism, although the latter construction always provides a canonical choice.\footnote{The relations of the form $x\cov_1 h([\{y\}])$ are those which satisfy the additional property $x\cov_1 r^{-1}\{y\}$ iff $r(x,y)$. In principle, this additional property could be required in the definition of a morphism; however, the composition of two such relations need not satisfy the same  property and one needs a somehow unnatural definition of composition. For this reason, we prefer to work with the usual composition of relations, though in doing so we must consider relations up to a suitable equivalence.}

If $(S_1,\cov_1)$ and $(S_2,\cov_2)$ are formal covers, a frame homomorphism from $\mathcal{P}(S_2)/=_{\cov_2}$ to $\mathcal{P}(S_1)/=_{\cov_1}$ corresponds to a relation $r$ which, in addition to the previous conditions, preserves finite meets, that is, $r^{-1}\{a\wedge b\}$ $=_{\cov_1}$ $\{x\wedge y\ |\ x\in r^{-1}\{a\}, y\in r^{-1}\{b\}\}$ and  $r^{-1}\{1\}$ $=_{\cov_1}$ $\{1\}$.

Now if $(S_1,\cov_1)$ and $(S_2,\cov_2)$ are $\sigma$-covers, frame homomorphisms which preserve Lindel\"of elements correspond to relations $r$ for which $r^{-1}\{b\}$ is countable ``up to the cover'' for all $b\in S_2$. Indeed, $h:\mathcal{P}(S_2)/=_{\cov_2}\to\P(S_1)/=_{\cov_1}$ preserves Lindel\"of elements if, and only if, for every $b\in S_2$, there is a countable $W\sub S_1$ with $h([\{b\}])=[W]$.

In view of the previous discussion, the following definition makes the category of $\sigma$-covers equivalent to the category of $\sigma$-locales.

\begin{defi}
A morphism between two $\sigma$-covers $(S_1,\cov_1)$ and $(S_2,\cov_2)$ is a binary relation $r\sub S_1\times S_2$ such that
\begin{enumerate}
\item if $U=_{\cov_2} V$, then $r^{-1}U=_{\cov_1}r^{-1}V$;
\item $r^{-1}\{1\}=_{\cov_1}\{1\}$;
\item $r^{-1}\{a\wedge b\}=_{\cov_1}\{x\wedge y\ |\ x\in r^{-1}\{a\}, y\in r^{-1}\{b\}\}$;
\item for every $b\in S_2$ there is some $W\in\P(S_1)$ such that $r^{-1}\{b\}=_{\cov_1} W$;
\end{enumerate}
and two such relations $r$ and $s$ are equivalent, that is, they are equal as morphisms if $r^{-1}\{b\}=_{\cov_1}s^{-1}\{b\}$ for all $b\in S_2$.
\end{defi}
One can check that $r^{-1}\{b\}=_{\cov_1}s^{-1}\{b\}$ for all $b\in S_2$ yields $r^{-1}U=_{\cov_1}s^{-1}U$ for all $U\sub S_2$. And the usual composition of relations is compatible with equality of morphisms. We note that showing that the composition
of morphisms satisfies property 4 requires $AC_\omega$.

\paragraph{Points.} The initial $\sigma$-frame $\Sigma$ is a terminal $\sigma$-locale; so the $\sigma$-locales arrows from it to a given $L$ are the (global) {\bf points} of $L$. In the case of locales, a point can be identified with a completely prime filters of opens. Similarly, a point of a $\sigma$-frame $L$ is a subset $p\sub L$ such that:
\begin{enumerate}
\item $1\in p$;
\item if $a\in p$ and $b\in p$, then $a\wedge b\in p$;
\item if $a\in p$ and $a\leq\bigvee W$ with $W\in\P(L)$, then $w\in p$ for some $w\in W$;
\item the truth value of $a\in p$ is in $\Sigma$.\footnote{More mathematically, this means that $\{x\in 1\ |\ a\in p\}$ is a countable subset of 1.}
\end{enumerate}

\noindent
When the category of $\sigma$-locales is embedded in the category of locales, as above, then a point of a $\sigma$-locale $L$ is the same thing as a point of its envelope $(L,\cov_L)$ under the proviso that $a\in p$ is ``semi-decidable'' for all $a\in L$.
This justifies the following definition.

\begin{defi}
A (global) {\bf point} of a $\sigma$-cover $(S,\cov)$ is a subset $p\sub S$ which satisfies the following
\begin{enumerate}
\item $1\in p$;
\item if $a\in p$ and $b\in p$, then $a\wedge b\in P$;
\item if $a\in p$ and $a\cov U$, then $u\in p$ for some $u\in U$;
\item the truth value of ``$a\in p$'' is in $\Sigma$.
\end{enumerate}
\end{defi}

\noindent
Note that a point of $(S,\cov)$ as a $\sigma$-cover is, in particular, a point of $(S,\cov)$ as a formal cover, but the converse fails, in general. And a morphism between $\sigma$-covers must map points (in this stronger sense) to points (in this stronger sense).

\subsection{Inductive generation of \texorpdfstring{$\sigma$}{sigma}-covers}

There are several important cases of formal covers which can be inductively generated. The general method is described in~\cite{cssv} (see also~\cite{convergence}). Although we are not going to give all details, the idea is to construct $\cov$ as the smallest cover which satisfies some set of ``axioms'' of the form $a\cov C(a,i)$, for $a\in S$ and $i$ in some given set $I(a)$.
When $S$ is assumed to have an inf-semilattice structure, as we always do in this paper, the cover generated by such a set of axioms is the smallest sub-collection $\cov\sub S\times\mathcal{P}(S)$ that satisfies the following clauses:
\begin{enumerate}
\item if $a\in U$, then $a\cov U$;
\item if $a\leq b$ and $b\cov U$, then $a\cov U$;
\item if $a\leq b$ and $c\wedge a\cov U$ for all $c\in C(b,i)$, then $a\cov U$;
\end{enumerate}
(we refer the reader to~\cite{cssv,convergence} for details).

\begin{prop} A cover $(S,\cov)$ is a $\sigma$-cover if and only if it can be inductively generated by means of axioms $C(a,i)$ which are all countable.
\end{prop}
\begin{proof}
Let $(S,\cov)$ be a $\sigma$-cover, put $I(a)=\{W\in\P(S)\ |\ a\cov W\}$ and $C(a,W)=W$ for $W\in I(a)$, and let $\cov'$ be the cover generated by these axioms. We claim that $\cov=\cov'$. Clearly $\cov'$ is contained in $\cov$ because $\cov'$ is the smallest cover satisfying the axioms, but $\cov$ satisfies them as well. Also, $\cov$ is contained in $\cov'$; for if $a\cov U$, then $a\cov W\subseteq U$ for some countable $W$; so $W\in I(a)$ and hence $a\cov'C(a,W)=W\subseteq U$.

Vice versa, let $(S,\cov)$ be inductively generated with all $C(a,i)$'s countable.
Assume $a\cov U$; we must show that $a\cov W$ for some countable $W\sub U$. The proof is by induction on the generation of the cover, of course; three cases can occur: (i) $a\in U$, (ii) there is $a\leq b$ with $b\cov U$, (iii) there is $b\geq a$ and $i\in I(b)$ with $c\wedge a\cov U$ for all $c\in C(b,i)$. The first case is trivial because $\{a\}$ is a countable subset of $U$. The second case is easy: we have $b\cov W$ for some countable $W\sub U$, by inductive hypothesis; so $a\cov W$ by clause 2. The third case requires $AC_\omega$: for each $c\in C(b,i)$ we have $c\wedge a\cov W_c\sub U$, with $W_c$ countable, by inductive hypothesis; let $W$ be the union of the $W_c$'s, which is a countable subset of $C(b,i)$ because $C(b,i)$ is countable by assumption; so $c\wedge a\cov W$ for all $c$, and hence $a\cov W$ by clause 3.
\end{proof}

Many important examples of generated covers are in fact $\sigma$-covers, such as (the point-free versions of) the Cantor space, the Baire space, and the (Dedekind) reals. Let us analyse the last example in details.

\paragraph{The reals.} The locale of the reals can be presented as follows. Let $\mathbb{Q}$ be the set of rational numbers. Put $S$ = $\{(a,b)\in\mathbb{Q}\times\mathbb{Q}\ |\ a<b\}$; this is a poset where $(a_1,b_1)\leq(a_2,b_2)$ means $a_2\leq a_1<b_1\leq b_2$. Note that $S$ is countable, because it is a detachable subset of a countable set (Lemma~\ref{lemmaDetachableCountable}). In order to turn $S$ into an inf-semilattice, we add a top element $(-\infty,+\infty)$; we also need a bottom element, say $(0,0)$, so that $(a_1,b_1)\wedge(a_2,b_2)$ is $(\max\{a_1,a_2\},\min\{b_1,b_2\})$ if $\max\{a_1,a_2\}<\min\{b_1,b_2\}$, and $(0,0)$ otherwise.

An element $(a,b)\in S$ is thought of as the open interval $]a,b[$. We want to define a cover $(a,b)\cov U$ in such a way to capture the intuition that $]a,b[\subseteq\bigcup_{(x,y)\in U}]x,y[$. This can be done by induction by means of the following axioms:
\begin{enumerate}
\item $(a,b)\cov\{(a,b'),(a',b)\}$ whenever $a<a'<b'<b$;
\item $(a,b)\cov\{(a',b')\ |\ a<a'<b'<b\}$.
\end{enumerate}
And the previous proposition applies, because $\{(a',b')\ |\ a<a'<b'<b\}$ is a detachable subset of the countable set $S$.
So this is a $\sigma$-cover and it makes sense to consider the $\sigma$-locale of its Lindel\"of elements, namely, the countable unions of basic opens. Classically, this is the whole frame of opens, because every open set can be written as a countable union of open intervals with rational endpoints;  constructively, since a subset of a countable set (the base, in this case) need not be countable, the $\sigma$-locale of Lindel\"of opens seems to keep its own independence.

As explained in the previous section, a point of this $\sigma$-locale is a particular filter $p$ of opens such that $(a,b)\in p$ is semi-decidable for every basic open $(a,b)$. So we can think of $p$ as a recursively enumerable set of basic opens: this is often called a computable point (see~\cite[Definition 2.2]{escardo} for example). Thus a $\sigma$-locale morphism between reals corresponds to the idea of a continuous function which maps computable points to computable points.

\section{On strongly dense \texorpdfstring{$\sigma$}{sigma}-sublocales}

A {\bf congruence} $\sim$ on a $\sigma$-frame $L$ is an equivalence relation which is compatible with finite meets and countable joins; this says that $L/\sim$ is a $\sigma$-frame as well. From the ``dual'' point of view of locale theory~\cite{stone_spaces}, $L/\sim$ is what is called a {\bf $\sigma$-sublocale} of $L$.\footnote{In the case of frames, congruences correspond to nuclei~\cite{stone_spaces}; such a correspondence relies on the existence of the implication operation (every frame has a Heyting algebra structure) and does not hold in the case of $\sigma$-frames. As in the case of frames, the family of congruences on a given $\sigma$-frame is a frame: we refer the reader to~\cite{madden,simpson} for more on this topic.}

A $\sigma$-sublocale $L/\sim$ is called {\bf dense} when $\forall x\in L.(x\sim 0\Rightarrow x=0)$. Every $\sigma$-locale $L$ has a smallest dense $\sigma$-sublocale\footnote{See~\cite[Proposition 6.2]{madden} where the $\sigma$-frames arising in this way are called \emph{d-reduced}.} which corresponds to the congruence $a\sim b$ defined by $\forall x\in L.(a\wedge x=0\Leftrightarrow b\wedge x=0)$.

By extending the notion given in~\cite{strong_density} from locales to $\sigma$-locales, we say that $L/\sim$ is a {\bf strongly dense} $\sigma$-sublocale of $L$ if
$x\sim\ !_L(p)$ implies $x\leq\ !_L(p)$, for all $x\in L$ and $p\in\Sigma$, where $!_L$ is the unique $\sigma$-frame homomorphism from $\Sigma$ to $L$. Clearly, strong density implies density by choosing $p=0$; and they coincide in a classical framework because $\Sigma=\{0,1\}$ in that case.

In what follows we are going to characterize the smallest strongly dense $\sigma$-sublocale of an \emph{overt} $\sigma$-locale. First we need to introduce the notion of overtness for $\sigma$-locales.

\subsection{Overt \texorpdfstring{$\sigma$}{sigma}-locales}

Recall that a formal cover $(S,\cov)$ is {\bf overt} if there is a predicate $\Pos(x)$ on $S$ (the {\bf positivity} predicate) such that
\begin{enumerate}
\item if $a\cov U$ and $\Pos(a)$, then $\Pos(b)$ for some $b\in U$,
\item $a\cov\{a\}\cap\Pos$,
\end{enumerate}
where $\Pos$ is $\{x\in S\ |\ \Pos(x)\}$. For instance, the locale of the reals is overt with $\Pos((a,b))$ whenever $a<b$ and (hence) also $\Pos((-\infty,+\infty))$.

The idea is that $\Pos(a)$ is a positive way to say that $a$ is not the bottom element. Note that $\neg \Pos(a)$ is (intuitionistically) equivalent to $a\cov\emptyset$ and hence also to $a=0$. Classically, therefore, $\Pos(a)$ means just $a\neq_\cov\emptyset$, and the conditions above are automatically satisfied.

It is well known, that a formal cover is overt precisely when the unique frame homomorphism $\Omega=\mathcal{P}(1)\to\mathcal{P}(S)/=_\cov$ has a left adjoint (see, for instance,~\cite{embedding}).

We want to extend this notion to $\sigma$-locales. We say that a $\sigma$-locale $L$ is {\bf overt} if its envelope $(L,\cov_L)$ is overt.\footnote{It is well known that for an overt $(S,\cov)$ the predicate $\Pos(a)$ becomes equivalent to the ``impredicative'' formula $\forall U\sub S.(a\cov U\Rightarrow U\textrm{ is inhabited})$. When $(S,\cov)$ is $\sigma$-coherent, that formula can be replaced by $\forall W\in\P(S).(a\cov W\Rightarrow W\textrm{ is inhabited})$.} Explicitly, $L$ is overt if and only if there is a predicate $\Pos(x)$ on $L$ such that
\begin{enumerate}
\item if $x\leq \bigvee W$ with $W$ countable and $\Pos(x)$, then $\Pos(w)$ for some $w\in W$;
\item for each $x$, there is a countable $W\sub\{x\}\cap\Pos$ such that $x\leq\bigvee W$.
\end{enumerate}
The second condition can be replaced with the simpler
\begin{enumerate}
\item[2'.] $\{x\}\cap\Pos$ is countable and $x\leq\bigvee(\{x\}\cap\Pos)$
\end{enumerate}
because if $W\sub\{x\}\cap\Pos$ and $x\leq\bigvee W$, then $W=\{x\}\cap\Pos$. To see this, let $y\in\{x\}\cap\Pos$; so $\Pos(x)$, and hence $\Pos(w)$ for some $w\in W$ by condition $1.$; in particular, $\{w\}\sub W\sub\{x\}$ and hence $w=x$; therefore $x\in W$. 

As a consequence of $2'$, we always have $\Pos(x)\in\Sigma$, because $\Pos(x)$ means that $\{x\}\cap\Pos$ is inhabited (Remark~\ref{remark_inhabited}). And, vice versa, if $\Pos(x)\in\Sigma$, then $\{x\}\cap\Pos$ is countable. Indeed $\Pos(x)$, seen as a subset of 1, is $\{z\in 1\ |\ \Pos(x)\}$; and there is a detachable $D\sub\N$ and a function $\alpha:D\to 1$ such that $\alpha[D]$ = $\{z\in 1\ |\ \Pos(x)\}$. If $k_x:1\to S$ is the constant map with value $x$, we have $k_x[\alpha[D]]$ = $k_x[\{z\in 1\ |\ \Pos(x)\}]$ = $\{y\in\{x\}\ |\ \Pos(x)\}$ = $\{y\in\{x\}\ |\ \Pos(y)\}$ = $\{x\}\cap\Pos$. Therefore we can replace $2'$ with
\begin{enumerate}
\item[2''.] $\Pos(x)\in\Sigma$ (hence $\{x\}\cap\Pos$ is countable) and $x\leq\bigvee(\{x\}\cap\Pos)$.
\end{enumerate}
It is routine to check that $\Pos$ is left adjoint to the unique $\sigma$-frame homomorphism $!_L:\Sigma\to L$, that is,
$\Pos(a)\Rightarrow p$ if and only if $a\leq\ !_L(p)$,
for all $a\in L$ and $p\in\Sigma$. And a $\sigma$-locale $L$ is overt if and only if $!_L$ has a left adjoint.

\begin{prop} Let $L$ be an overt $\sigma$-locale and let $\sim$ be a congruence on $L$.
Then $L/\sim$ is strongly dense if and only if $\Pos(x)\Rightarrow\Pos(y)$ whenever $x\sim y$.
\end{prop}
\begin{proof}
Let $L/\sim$ be strongly dense and assume $x\sim y$. Since $y\leq\ !_L(\Pos(y))$, we also have $(x\vee\ !_L(\Pos(y)))\sim\ !_L(\Pos(y))$. So $x\vee\ !_L(\Pos(y))\leq\ !_L(\Pos(y))$ and hence $x\leq\ !_L(\Pos(y))$, that is, $\Pos(x)\Rightarrow\Pos(y)$.

Vice versa, if $x\sim\ !_L(p)$, then $\Pos(x)\Rightarrow\Pos(!_L(p))$ by assumption; now the usual properties of an adjunction gives $x\leq\ !_L(\Pos(!_L(p)))=!_L(p)$.
\end{proof}

\subsection{Overlap algebras}

Complete Boolean algebras lose some of their important features when LEM is not assumed. For instance, discrete locales, that is, frames of the form $\mathcal{P}(S)$ for some set $S$, are never Boolean, apart from the trivial case $S=\emptyset$.\\
Sambin's notion of an \emph{overlap algebra} (see below) is a constructive  alternative to that of a complete Boolean algebra (see~\cite{regolari, formal_regular} for some basic results). For instance, powersets are examples of overlap algebras, in fact they are precisely the atomic ones.\\
It has recently turned out~\cite{almost_discrete} (see also~\cite{contente}) that overlap algebras can be understood as the smallest strongly dense sublocales of overt locales. Under a classical reading, this means that overlap algebras are precisely the Boolean locales (since these can be
characterized as the smallest
dense sublocales).\\ Such a result can be adapted to characterize the strongly dense $\sigma$-sublocales, which is the main point of the current
section.

An {\bf overlap algebra} is (the locale corresponding to) an overt formal cover $(S,\cov,\Pos)$ such that
\begin{equation}\label{Posa}
\forall b\in S.[\Pos(a\wedge b)\Rightarrow\exists u\in U.\Pos(u\wedge b)]\Longrightarrow a\cov U
\end{equation}
for every $a\in S$ and $U\sub S$.\footnote{Actually, this is the definition of a \emph{set-based} overlap algebra: for the general definition one has to consider an arbitrary locale, not necessarily corresponding to a formal cover (such a distinction makes sense in a predicative framework).}

\begin{defi}
A {\bf $\sigma$-overlap algebra} is an overt $\sigma$-locale $L$ such that
\begin{equation}\label{oa}
\forall b\in L.[\Pos(a\wedge b)\Rightarrow \Pos(u\wedge b)]\Longrightarrow a\leq u
\end{equation} for all $a,u\in L$.
\end{defi}
Classically, $\sigma$-overlap algebras are precisely the d-reduced $\sigma$-frames in the sense of~\cite[Section 6]{madden} where they play an important role in connection with \emph{regularity} (see~\cite[Proposition 6.5]{madden} for instance).

If $S$ is countable with semi-decidable equality, then $\P(S)$ is a $\sigma$-overlap algebra where $\Pos(W)$ is ``$W$ is inhabited'' (see Propositions~\ref{prop1} and~\ref{prop2} and Remark~\ref{remark_inhabited}). More generally, if the envelope $(L,\cov_L)$ of a $\sigma$-frame $L$ is an overlap algebra, then $L$ is a $\sigma$-overlap algebra; indeed,~\eqref{oa} is clearly a special case of~\eqref{Posa} in that case.

Further examples of $\sigma$-overlap algebras can be constructed as follows.
Let $L$ be an overt $\sigma$-locale and let $\sim_{B(L)}$ be the binary relation on $L$ defined by
\begin{equation}\label{o}
x\sim_{B(L)} y\ \stackrel{def}{\Longleftrightarrow}\ \forall z\in L.(\Pos(x\wedge z) \Leftrightarrow \Pos(y\wedge z)).
\end{equation}
It is quite straightforward to check that~\eqref{o} defines a congruence on $L$, hence a $\sigma$-sublocale of $L$.

\begin{defi}
Given any overt $\sigma$-locale $L$, we write $B(L)$ for the $\sigma$-sublocale of $L$ corresponding to the congruence $\sim_{B(L)}$ as defined in~\eqref{o}.
\end{defi}
Note that $[x]\leq[y]$ in $B(L)$ if and only if $\Pos(x\wedge z)\Rightarrow\Pos(y\wedge z)$ for all $z\in L$. Indeed $[x]\leq[y]$ iff $x$ $\sim_{B(L)}$ $x\wedge y$ iff $\Pos(x\wedge z) \Leftrightarrow \Pos(x\wedge y\wedge z)$ for all $z$
 iff $\Pos(x\wedge z) \Rightarrow \Pos(x\wedge y\wedge z)$ for all $z$ iff $\Pos(x\wedge z) \Rightarrow \Pos(y\wedge z)$ for all $z$.\footnote{The last step holds because if $\Pos(x\wedge z) \Rightarrow \Pos(y\wedge z)$ for all $z$, then also $\Pos(x\wedge x\wedge z) \Rightarrow \Pos(y\wedge x\wedge z)$ for all $z$.}

The $\sigma$-sublocale $B(L)$ is always dense in $L$. Indeed, if $a\sim_{B(L)}0$, that is, $\forall z[\Pos(a\wedge z) \Leftrightarrow \Pos(0\wedge z)]$, then $\forall z.\neg \Pos(a\wedge z)$ because $\Pos(0)=\Pos(\bigvee\emptyset)$ is false. In particular, $\neg \Pos(a)$ and so $a=0$. Actually, we have the following.

\begin{prop}
If $L$ is an overt $\sigma$-locale, then
\begin{enumerate}
\item $B(L)$ is the smallest strongly-dense $\sigma$-sublocale of $L$;
\item $B(L)$ is a $\sigma$-overlap algebra.
\end{enumerate}
\end{prop}
\begin{proof}\hfill
\begin{enumerate}
\item Let $L/\sim$ be strongly dense. If $x\sim y$, then $(x\wedge z)\sim(y\wedge z)$ because $\sim$ is a congruence; hence $\Pos(x\wedge z)\Leftrightarrow \Pos(y\wedge z)$ by strong density. This means that $B(L)$ is a $\sigma$-sublocale of $L/\sim$.
\item First we check that $B(L)$ is overt. This is easy because $\Pos$ respects the congruence~\eqref{o} and so it makes sense to define $[a]$ positive in $B(L)$ if $a$ is positive in $L$. The two conditions on $\Pos$ are easy to check.\\
Now, if $x$ and $y$ are such that $\Pos(x\wedge z) \Rightarrow \Pos(y\wedge z)$ holds in $B(L)$ for all $z$, then $x\leq y$ in $B(L)$.\qedhere
\end{enumerate}
\end{proof}

\noindent When $L$ is a $\sigma$-overlap algebra, we have $B(L)=L$ because $\sim_{B(L)}$ becomes the identity in that case. So we immediately have the following.

\begin{cor} A(n overt) $\sigma$-locale $L$ is a $\sigma$-overlap algebra if and only if $L$ is (isomorphic to) the smallest strongly-dense $\sigma$-sublocale $B(X)$ of some overt $\sigma$-locale $X$.\qed
\end{cor}

\subsection*{Acknowledgements}
The author wishes to thank two anonymous referees for their interesting comments, their useful suggestions, and their challenging questions.

\bibliographystyle{alphaurl}
\bibliography{ciraulo}

\end{document}